\documentstyle[11pt]{article}\topmargin=-2cm\textheight=210mm
\begin{document}\bibliographystyle{unsrt}
\begin{center}{\LARGE \bf Classical Mechanics of Spinning Patricle\\
in a Curved Space}\\[5mm] Z.\,Ya.\,Turakulov\\{\it Institute of Nuclear
Physics\\ Ulugbek,Tashkent 702132, Rep. of Uzbekistan, CIS\\
(e-mail: zafar@suninp.tashkent.su)}\\[5mm] \end{center}\begin{abstract}
An example of mechanical system whose configuration space is direct
product of a curved space and the local group of rotations, is presented.
The system is considered as a model of spinning particle moving in the
space. The Hamiltonian formalism for this system and possible method for
its quantization are discussed. It is shown that the Hamilton equations 
coincide with the Papapetrou equations for spinning test-particle in general 
relativity.\end{abstract}
\section{Introduction}

Motion of spinning objects in curved spaces has been studied by numerous
authors \cite{{pu},{scf},{sff},{bmss}} as a purely general relativistic
problem. However, one special case of this phenomenon has been described 
exhaustively in standard texts on classical mechanics. In order to show 
this consider the simplest case of such a motion assuming that the space
is a two-dimensional sphere and the moving object is a two-dimensional
rigid body covering a domain in this sphere. As the sphere is a symmetric
space and, hence, it admits free motions of a rigid body in it, the body
represents a kind of spinning objects moving in this space. On the other
hand, this mechanical system imbedded into three-dimensional Euclidean
space constitutes a case of three-dimensional rigid body with fixed
point. Indeed, the center of the sphere can be considered as a point
the body because all the distances between this and other points of the
body remain constant thus, inclusion of the center does not break the
rigidity condition.  Consequently, the well-known branch of mechanics
studying motions of rigid bodies constitutes a primer of mechanics of
rotating objects in curved spaces.  In the present work we compose
non-relativistic Hamiltonian description of mechanical systems formed
by spinning objects moving in an arbitrary curved space.

The first problem arising when considering motions of a rigid body in an
arbitrary curved space is that, in general, such a space does not admit
free motions of a rigid body covering some domain on it. In order to
eliminate this difficulty we suggest some generalization of the conventional
definition of rigid body. Due to our definition a rigid body with center of
mass in some point $x$ of the space covers a domain in the tangent
space $T_x$ over this point \cite{klg} and therefore we will call it
``tangent rigid body''. Apparently, such a body is able to move freely on
a curved space regardless of intrinsic geometry of this space.

Below we will distinguish two main kinds of motion of tangent rigid body.
One of them is motion of the center of mass under which the body obeys
the law of parralel transfer. Another one represents pure rotation, i.e.,
the motion under which the center of mass rests. As was pointed out in
the work \cite{bmss} the configuration space of the system is locally
diffeomorfic to the group of symmetry of the tangent space (the Poincare
group for the case considered in the cited work). Therefore, its subspaces
containing trajectories of the main motions, are orthogonal. The
orthogonality means that the kinetic energy of the system is equal to the
sum of kinetic energies of the main motions. Nevertheless, these orthogonal
motions depend on each other and their interaction yields a difference
between behaviours of tangent rigid bodies and mass points.

In order to demonstrate this interaction consider the following special case
of tangent rigid body. Let the space will be specified as a surface in
three-dimensional Euclidean space and the body is a plate remaining tangent
to the surface. Then, the vector of angular momentum of the body, concerned
with its rotation, lies on the normal to the surface established in the center
of mass of the body. Since any infinitesimal displacement of the body yields an
infinitesimal rotation of the vector of angular momentum it yields a moment of
force acting on the body. In the simplest cases the moment of force is
collinear with the displacement. The constraints confining the body in tangent
planes transform this moment into a force damping body tranversally to the
original displacement. Consequently, the force acting on a body depends
linearly on its velocity, angular momentum and the space curature, whereas its 
direction is orthogonal to the velocity. Below we generalize these 
consideratons to the case of arbitrary Riemann space.\section{Definition of
tangent rigid body}

Consider a Riemann space $M$ endowed with a coordinate system $\{\,x^i\}$ and 
let 1-forms $\theta^a$ \begin{equation} \theta^a=e_i^a(x)dx^i \end{equation} 
constitute an orthonormal frame in every point of the space $M$ such that
its metric tensor may be expressed through the Kronecker delta:
$$g=g_{ij}dx^i \otimes dx^j= \delta_{ab} \theta^a \otimes \theta^b.$$
The connection 1-form for these frames may be introduced through the first
structure equation:\begin{equation}d\theta^a=\omega_b^a\wedge\theta^b,
\quad \omega_{ab}+\omega_{ba}=0\end{equation}and the connection
coefficients $\gamma_{abc}$ are that of expansion of this 1-form in the local 
frames $\{\,\theta^{a}\}$:\begin{equation} \omega_{b.}^c=\gamma_{ab.}^c 
\theta^a.\end{equation} The Riemann curvature is represented by the 2-form 
$\Omega_{a.}^{b}$ defined by the second structure equation \begin{equation} 
\Omega_{a.}^b=d \omega_{a.}^b+ \omega_{a.}^c \wedge \omega_{c.}^b
\end{equation}and the coefficients of its expansion in the local frame of
2-forms $\{\,\theta^{a} \wedge \theta^{b} \}$ constitute the components of
the Riemann tensor of curvature:\begin{equation}
\Omega_{a.}^b=R_{cda.}^b\theta^c \wedge \theta^d.\end{equation}

By analogy with classical definition of rigid body \cite{ad} the notion
of tangent rigid body may be defined as a collection of mass points
moving in the tangent space under holonomy constraints due to which
the distance between any two of them remains constant. Thus, let $T_x$
be a tangent space over a point $x$, $m_{\alpha}$ and $r_{\alpha}\in T_x$
be masses and radius-vectors of the points forming the body. Then,
apparently, the mass of the body is equal to\begin{equation}
m=\sum_{\alpha} m_{\alpha}\end{equation}Assuming that the center of mass
of the body coincides with the point $x\in M$ one obtains the following
identity:\begin{equation}\sum_{\alpha} m_{\alpha}r_{\alpha}=0.
\end{equation}The following sum represents components of the tensor of
inertia of the body:\begin{equation}
I^{ab}=\sum_{\alpha}m_{\alpha}r^a_{\alpha}r^b_{\alpha}\end{equation}
and it will be assumed below that the body is symmetric such that this
tensor has the form of Kronecker delta:\begin{equation}I^{ab}=I\delta^{ab}.
\end{equation}

As well as in the classical case \cite{ad} the configuration space of the
tangent rigid body is the direct product of space $M$ and the group of
rotations of local frames. Coordinates of the system in this space are
specified by coordinates of the center of mass $x^i$ and by orientation of
the body with respect to the local frame $\{\,\theta^a\}$. Since none of
forces act in this system its Lagrangian and Hamiltonian are equal to its
kinetic energy \cite{ad}. Now, in order to compose Lagrangian and
Hamiltonian formalisms for this system we introduce generalized velocities 
concerned with coordinates of the system in the configuration space.
\section{Generalized coordinares and generalized velocities}

As was pointed out above the generalized coordinates of the system in question
are coordinates $x^i$ of the body center of mass and some functions of
components of the vectors $r_{\alpha}$ referred to the local frame.
Consequently, the corresponding generalized velocities are $\dot x^i$ and
some functions of components of the vectors $r_{\alpha}$ and their time
derivatives ${d\over dt}r_{\alpha}^a$, depending on the derivatives
linearly. The rigidity conditions require existance of antisymmetric
tensor $\eta_{a.}^b$ of angular velocity specifying all the time
derivatives of $r_{\alpha}^a$'s as follows:
$${d\over dt}r_{\alpha}^a=\eta_{b.}^ar_{\alpha}^b.$$The components of this
tensor may be considered as generalized velocities of the system. Since the
kinetic energy of the body does not depend on its orientation the
corresponding generalized coordinates may be omitted such that the
Lagrangian is a function of variables $x^i$, $\dot x^i$ and $\eta_{a.}^b$
only. It must be pointed out that definition ofangular velocities introduced
above do not account changes of local frames along the center of mass
trajectory. Therefore, it is necessary to introduce the covariant time
derivative of $r_{\alpha}$ components:\begin{equation}
\dot r_{\alpha}^a\equiv(\eta_{b.}^a+\dot x^c\gamma_{cb.}^a)r_{\alpha}^b.
\end{equation}
In order to compose the expression of Lagrangian consider an infinitesimally
small tangent rigid body assuming that all the mass points forming the body
lie in a small neighbourhood of the center of mass and have coordinates
$x+r_{\alpha}$. Then their velocities are $\dot x+ \dot r_{\alpha}$ with
$\dot r_{\alpha}$ being defined in the equality (10). Now, taking into 
account the equations (6-10) one obtains the following expression for kinetic 
energy:  \begin{equation} L=T={1\over 2}\sum_{\alpha}m_{\alpha}(\dot x+\dot 
r_{\alpha}, \dot x+\dot r_{\alpha})= {m\over 2}(\dot x,\dot x)+{1\over 
2}I(\eta_{a.}^a+\gamma_{cb.}^b\dot x^c) (\eta_{a.}^b+\gamma_{c.a}^b\dot x^c) 
\end{equation}This expression conveys only the first non-vanishing
approximation for the kinetic energy. In this case the further analysis is
related to the system defined by the obtained form of Lagrangian. It must
be pointed out that non-orthogonality of this form does not contradict the 
above statement that the main kinds of motion are orthogonal, because, by 
definition, the angular velocities $\eta_{a.}^b$ are non-covariant.
\section{Generalized momenta and Poisson brackets}

The direct evaluation of the generalized momenta from the expression(11)
yields\begin{eqnarray}
p_i={\partial L\over{\partial\dot x^i}}=[(m\delta_{ab}+I\gamma_{ac.}^d
\gamma_{a.d}^c)\dot x^b+I\gamma_{abc}\eta^{bc}]e_i^a;\quad
S_{ab}={\partial L \over\partial\eta^{ab}}=I(\eta_{ab}+\gamma_{cab}
\dot x^c)\end{eqnarray}Here $e_i^a$ is the matrix introduced in the equation
(1), $p_a$ are components of the body momentum and $S_{ab}$ are that of its 
spin. The Poisson brackets for the spin components coincide with commutators
of corresponding generators of the group of rotations:\begin{equation}
[S_{ab},S_{cd}]={1\over2}(\delta_{ac}S_{bd}+\delta_{bd}S_{ac}-\delta_{ad}
S_{bc}-\delta_{bc}S_{ad})\end{equation}Other Poisson brackets for
generalized momenta are, by definition, equal to zero:\begin{equation}
[p_i,p_j]=[p_i,S_{ab}]=0\end{equation}whereas $[p_a,p_b]$ is not zero
since the orthonormal components $p_a$ contain the factors $e_a^i$ depending
on $x^i$'s.

In order to obtain the explicit form of the Hamiltonian one must express
generalized velocities through the generalized momenta. The direct
evaluation yields:$$I\eta_{ab}=S_{ab}-\gamma_{cab}\dot x^c.$$
Substituting this into the first equality (12) one obtains:\begin{equation}
m\dot x^a=\delta^{ab}(p_b-\gamma_{bcd}S^{cd}).\end{equation}The last two
results transform the expression (11) to the following form:\begin{equation}
H={1\over2m}\delta^{ab}(p_a-\gamma_{acd}S^{cd})(p_b-\gamma_{bcd}S^{cd})
+{S^2\over2I}\end{equation}\section{Equations of motion}

As follows from the second equation (12) and independence of the 
Lagrangian on the body orientation the spin $S_{ab}$ remains constant
on tha trajectory. Consequently, the last term in the equation (16) is
constant and, hence, it can be substracted from the Hamiltonian, i.e.,
the Hamiltonian can be written in the form\begin{equation}
H={1\over2m}\delta^{ab}(p_a-\gamma_{acd}S^{cd})(p_b-\gamma_{bcd}S^{cd})
\end{equation}It is to be pointed out that the equality (15) is
conspicuously similar to the well-known definition of the covariant
momentum of a particle moving in a gauge field and the expression (17)
coincides with that of Hamiltonian of such a particle. Remarkable, the
field in this analogy is represented by connection coefficients, whereas
the particle spin represents its charge coupled to the field.

Derivations of the Hamilton equation for spin leads, with account of the
Poisson brackets (13) and the equation (15), to the following equality:
$${dS_{ab} \over dt}=[H,S_{ab}]=-\dot x^c \gamma_c^{de}[S_{de},S_{ab}]=
-\dot x^c(\gamma_{ca.}^d S_{db}+\gamma_{cb.}^d S_{ad}),$$or, in agreement
with the Papapetrou's result,\begin{equation}{DS_{ab} \over dt}=0
\end{equation}Now, it remains to analyze the Hamilton equation
\begin{equation}{dp_a\over dt}=[H,p_a]\end{equation}It must be pointed out
that in the case of zero spin this equation coincides with that of geodesic
and, consequently, all its spin-independent terms constitute the covariant
time derivative of the vector $m \dot x$. The left hand side of this equation
is$${dp_a \over dt}=m\delta_{ab}{d\dot x^b \over dt}+\dot x^d \partial _d
\gamma_{abc}S^{bc}+\gamma_{abc}{dS^{bc} \over dt}$$where the derivatives
$\partial_a$ are defined as follows:$$\partial_a=e_a^i{\partial\over\partial
x^i}.$$Substituting here the time derivative of $S^{bc}$ from thr equation
(18) one finds that$${dp_a\over dt}=m\delta_{ab}{d\dot x^b\over dt}+\dot
x^d\gamma_{abc}(\gamma_{d.e}^b S^{ec}+\gamma_{d.e}^c S^{be}).$$The right-hand
side of the equation (19) contains the only spin-depending term$$\dot
x^d[p_a,\gamma_{dbc}]S^{bc}=(\partial_a\gamma_{dbc})\dot x^dS^{bc}.$$Now,
collecting all the evaluated terms and taking into account the fact that
all the spin-independent terms yield the covariant time derivative of
$m \dot x$ and the definitions of connection and curvature (2-5) one finds
that the Hamilton equation (19) coincides with the Papapetrou equation
\cite{pu}$$m{D \dot x^a\over dt}=R_{.bcd}^a\dot x^bS^{cd}.$$\section{Remarks
on quantization}

Quantization of the mechanical system considered above changes it into
a quantum mechanical description of motion of a particle with spin in
the space $M$. Operators of the spin components $\hat S_{ab}$ form a
representation of the Lie algebra of the group of rotations and, hence,
the wave function takes its values in the vector space of this
representation. Now, chosing the operators of the particle momentum
componenets in accord with the Poisson brackets (14)$$p_i=
{1\over\imath}{\partial\over\partial x^i}$$one finds that the Hamiltonian
operator contains only covariant derivatives of the wave function:$$\hat
H=g^{ij}({1\over\imath}{\partial\over\partial x^i}+e_i^a\gamma_{abc}\hat
S^{bc})({1\over\imath}{\partial\over\partial x^j}+e_j^a\gamma_{abc}\hat
S^{bc}).$$In other words, the particle momentum is represented by the
operator of covariant derivative$$\hat P_i={1\over\imath}{\partial\over
\partial x^i}+e_i^a\gamma_{abc}\hat S^{bc}$$instead of the operator of
generalized momentum introduced above.
\end{document}